\documentclass[12pt]{amsart}

\usepackage{graphicx}
\usepackage{amsfonts}
\usepackage{epsf}
\usepackage{amssymb}
\usepackage{amsmath}
\usepackage{amscd}
\usepackage{hyperref}
\usepackage{color}
\usepackage{bbm}
\usepackage{caption}
\usepackage{subcaption}

\newtheorem{theorem}{Theorem}[section]
\newtheorem*{utheorem}{Theorem}
\newtheorem{proposition}[theorem]{Proposition}

\newtheorem{corollary}[theorem]{Corollary}

\theoremstyle{remark}

\newtheorem{example}[theorem]{Example}
\newtheorem{remark}{Remark}[section]

\newcommand{\cokerr}{\mathrm{Coker}}

\input xy
\xyoption{all}

\makeatletter
\@namedef{subjclassname@2020}{\textup{2020} Mathematics Subject Classification}
\makeatother

\begin{document}

\title{Band-unknotting numbers and connected sums of knots}

\begin{abstract}
We study the band-unknotting number $u_{nb}(K)$ of a knot $K$, and how it behaves with respect to connect sums. We show that this sub-additive function is not additive under connected sums, by finding infinitely many examples of knots $K_1, K_2$ with $u_{nb}(K_1\#K_2) < u_{nb}(K_1) + u_{nb}(K_2)$. Even more surprisingly, there are infinitely many examples of knots $K_1, K_2$ such that $u_{nb}(K_1\#K_2) < u_{nb}(K_i)$, $i=1,2$. Our work is motivated by the recent analogous results for the Gordian unknotting number by  Brittenham and Hermiller \cite{BrittenhamHermiller}. 

We also prove new lower and upper bounds on the topological and smooth non-orientable 4-genus of a knot $K$. 
\end{abstract}

\author{Nakisa Ghanbarian}
\email{nakisa@unr.edu}
\author{Stanislav Jabuka}
\email{jabuka@unr.edu}
\address{Department of Mathematics and Statistics, University of Nevada, Reno NV 89557, USA.}
\maketitle
\section{Introduction}
\subsection{Background} An {\em unknotting operation $\mathcal O$} 
is a ``local'' transformation on a knot diagram yiedling another knot diagram, such that any diagram can be unknotted after finitely many steps.  The {\em $\mathcal O$-unknotting number $u_{\mathcal O}(K)$} is the minimum number of applications of the operation $\mathcal O$, taken across all diagrams of $K$, that yields a diagram of the unknot. It is easy to see that for oriented knots $K_1$ and $K_2$, the $\mathcal O$-unknotting number is subadditive with respect to the connected sum operation, that is
\begin{equation} \label{Universal inequality for unknotting numbers under connected sum}
u_{\mathcal O}(K_1\#K_2) \le u_{\mathcal O}(K_1) + u_{\mathcal O}(K_2).	
\end{equation}
It is an interesting question to determine if strict inequality can occur.

\begin{figure}[ht]
	\includegraphics[width=0.80\textwidth]{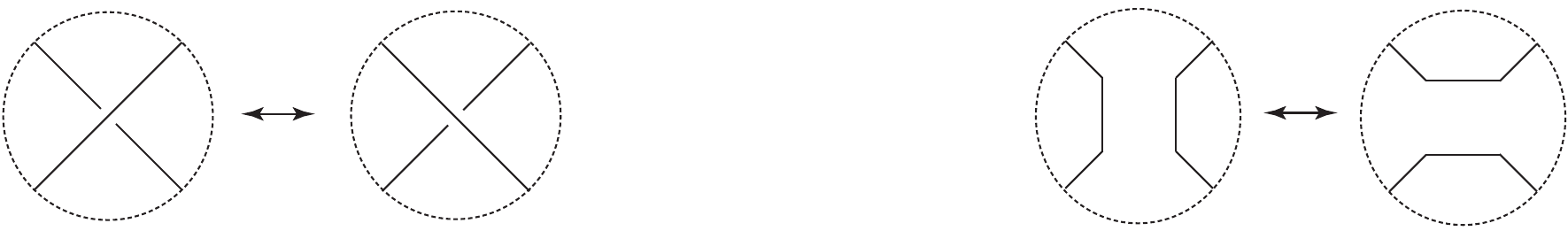}
	\put(-297,-15){(a)}
	\put(-67, -15){(b)}
	\caption{(a) The unknotting operation of switching a crossing. (b) The unknotting operation of a non-orientable band move. We require both diagrams to represent a knot. } \label{Unknotting operations}
\end{figure}

Classically, the most studied unknotting operation is that of a crossing-changing (Figure \ref{Unknotting operations}(a)) leading to the {\em Gordian or classical unknotting number $u(K)$}. Many other unknotting operations exist, see for example \cite{HosteNakanishiTaniyama} or Chapter 11 of \cite{Kawauchi}. 

A beautiful and remarkable recent result due to Brittenham and Hermiller \cite{BrittenhamHermiller} proves that inequality \eqref{Universal inequality for unknotting numbers under connected sum} can become sharp for the Gordian unknotting number.

\begin{utheorem}[Theorem 1.2, \cite{BrittenhamHermiller}]
Unknotting number is not additive under connected sum. In particular, for the knot $K = 7_1$ (the (2,7)-torus knot), with
$u(K) = 3$, the connected sum $L = K\#\bar K$ satisfies $u(L) \le 5 < 6 = u(K) + u(\bar K)$.  
\end{utheorem}
In the above, and throughout this work, $\bar K$ represents the mirror knot of $K$. Similarly, if $K$ is oriented, we let $-K$ represent the reverse of $K$. 
\vskip3mm
\subsection{The band-unknotting number} \label{Subsection from the introduction on band unknotting numbers}
In this work we focus on the unknotting operation given by performing {\em non-orientable band moves} on knot diagrams, as defined by Figure \ref{Unknotting operations}(b). It is easy to see that non-orientable band moves are an unknotting operation as a crossing-change can be obtained by two non-orientable band moves. We denote the corresponding unknotting number by $u_{nb}$, and refer to it as the {\em band-unknotting number}. Our main results are inspired by the work of Brittenham and Hermiller \cite{BrittenhamHermiller}.
\begin{theorem} \label{First theorem on inequality among band-unknotting numbers}
There exist infinitely many oriented knots $K_1$ and $K_2$ such that 
$$ u_{nb}(K_1\# K_2) < u_{nb}(K_1) + u_{nb}(K_2).$$
\end{theorem}
\noindent Example \ref{Example of knots for strict inequality for band unknotting numbers of connected sums} showcases some concrete knots realizing this strict inequality. 
\vskip1mm
In \cite{BrittenhamHermiller} the question is asked (Question 4.2) whether it is true for all knots $K_1, K_2$ that
$$u(K_1\#K_2)\ge \max\{u(K_1), u(K_2)\}.$$ 
This question remains unanswered for the Gordian unknotting number, but we are able to offer this answer for the analogous question for the band-unknotting number. 
\begin{theorem} \label{Theorem about strict inequality for max band unknotting number}
There exists infinitely many oriented knots $K_1, K_2$ such that 
$$u_{nb}(K_1\#K_2) < u_{nb}(K_i)\, \text{ for } \,  i=1,2.$$
\end{theorem}
The knot group is a powerful knot invariant. While not a complete knot invariant, it is known that prime knots with isomorphic groups are equivalent \cite[Corollary 2.1]{GordonLuecke}. It was shown in \cite[Corollary 1.4]{BrittenhamHermiller} that the knot group does not in general determine the Gordian unknotting number. We are able to prove the analogous result for the band-unknotting number. 
\begin{theorem} \label{Theorem about how the knot group does not determine the band unknotting number}
There exist infinitely many alternating, oriented knots $K_1, K_2$ such that 
$$u_{nb}(K_1\#K_2) =1 \quad \mathrm{ and } \quad u_{nb}(K_1\#(-\bar K_2))\ge 2.$$ 
Note that $\pi(K_1\#K_2) \cong \pi(K_1\#(-\bar K_2))$. One such example is given by $K_1=5_2$ and $K_2=6_1$, both of which are twist knots. 
\end{theorem}
\subsection{On non-orientable 4-genera} 
The results in Section \ref{Subsection from the introduction on band unknotting numbers} are implied by existing upper and lower bounds on $u_{nb}$ (reviewed in Section \ref{Section on background on band-unknotting numbers}) and new bounds presented in this section. To formulate them, we need these preliminaries. 
\vskip1mm
The {\em fusion number} $F(K)$ of a ribbon knot $K$, defined in \cite{JuhaszMillerZemke}, is the minimum integer $n\ge 0$, taken across all ribbon disks in $S^3$ bounded by $K$, which become the $(n+1)$-component unlink after attaching $n$ oriented bands to $K$. Attaching an oriented band is the same operation as described by Figure \ref{Unknotting operations}(b), with the caveat that when applied to a single component in a link, it increases the number of components by 1. It is shown in \cite[Section 1.6]{JuhaszMillerZemke} that
\begin{equation} \label{Upper bound on the fusion number in terms of the bride index}
F\left(K\#(-\bar K)\right) \le \textrm{br}(K)-1,
\end{equation}   
with br$(K)$ being the bridge index of $K$. 
\vskip1mm
The {\em topological} and {\em smooth non-orientable 4-genera} $\gamma_{4,t}(K)$ and $\gamma_{4,s}(K)$ of a knot $K$, are the minimum of the first Betti numbers $b_1(\sigma)$ with $\sigma$ ranging over all non-orientable, compact surfaces $\sigma$, properly embedded in the 4-ball $D^4$ with $\partial \sigma = K$. For $\gamma_{4,t}$ we require $\sigma$ to possess a tubular neighborhood (or equivalently, to be locally topologically flatly embedded) while for $\gamma_{4,s}$ we ask $\sigma$ to be smoothly embedded. Clearly $\gamma_{4,t}(K) \le \gamma_{4,s}(K)$ for every knot $K$. 
\vskip1mm 
For a rational homology 3-sphere $Y$, consider the linking form $\lambda :H_1(Y)\times H_1(Y)\to \mathbb Q/\mathbb Z$. The Witt Decomposition Theorem (see Section \ref{Section on background on linking forms} for more details) says that there is a decomposition of the form $(H_1(Y), \lambda)$ into an orthogonal direct sum $(H_1(Y), \lambda) \cong (G_1, \lambda_1)\oplus (G_2,\lambda_2)$ with $(G_1, \lambda_1)$ an anisotropic form (a form for which $\lambda_1(g,g)\ne 0$ for any nonzero $g\in G_1$), and with $(G_2,\lambda_2)$ metabolic. In this decomposition, the isomorphism type of $(G_1, \lambda_1)$ is determined by $(H_1(Y),\lambda)$. For the case of $Y=\Sigma_2(K)$ - the two-fold cyclic cover of $S^3$ branched along $K$ - let $\mu_{an}(K)$ denote the minimal number of generators for $G_1$. 
\begin{theorem} \label{Main double inequality for the band-unknotting number}
Let $K$ be an oriented knot. Then:
\begin{itemize}
	\item[(a)] $\mu_{an}(K)\le \gamma_{4,t}\le \gamma_{4,s}(K) \le u_{nb}(K)$.
	\item[(b)] $u_{nb}(K\#(-\bar K)) \le 2 F\left(K\#(-\bar K)\right)$. 
\end{itemize}
\end{theorem}
\begin{example} \label{Example of knots for strict inequality for band unknotting numbers of connected sums}
	According to \cite{KnotInfo}, any knot $K$ from among the list
	\begin{align*}
		4_1,\, 6_3, \, 7_5,\, 7_7, \, 8_1, \, 8_2, \, 8_{12}, \, 8_13,\, 9_2, \, 9_{10},\, 9_{11},\, 9_{12},\, 9_{14}, \, 9_{18}, \, 9_{20}, \, 10_2, \, 10_5, \, 10_9, \cr
		10_{10}, \, 10_{13}, \, 10_{14}, \, 10_{18}, \, 10_{19}, \,  10_{25},\, 10_{26}, \, 10_{28}, \, 10_{32}, \, 10_{33}, \, 10_{34}, \, 10_{36}, \, 10_{37},   
	\end{align*} 
	has bridge index and $\gamma_{4,s}$ equal to 2 (the values of $\gamma_{4,s}$ for 8--10 crossing knots were computed in \cite{Ghanbarian, JabukaKelly}). By Theorem \ref{Main double inequality for the band-unknotting number} and \eqref{Upper bound on the fusion number in terms of the bride index} we find
	\begin{align*}
		u_{nb}(K)= u_{nb}(-\bar K) & \ge \gamma_{4,s}(K) = 2, \qquad \text{ and }  \cr 
		u_{nb}(K\#(-\bar K)\le 2F\left(K\#(-\bar K)\right) & \le 2(\textrm{br}\left(K\#(-\bar K \right)-1) = 2.
	\end{align*} 
	We tacitly relied on the fact that the bridge number minus 1, is additive under connected sum, a classical results of Schubert's \cite{Schubert}. Taking $K_1=K$ and $K_2=-\bar K$ for any knot $K$ above, gives an example for Theorem \ref{First theorem on inequality among band-unknotting numbers}, since 
	$$u_{nb}(K_1\#K_2) \le  2 < 4 = 2+2 \le u_{nb}(K_1) + u_{nb}(K_2).$$
\end{example}
\begin{corollary} \label{Corollary on the upper bound of gamma4 by the unknotting number plus one}
For any knot $K$, the following inequalities holds: 
$$\mu_{an}(K)\le \gamma_{4,t}\le \gamma_{4,s}(K)\le u(K)+1.$$ 
\end{corollary}
These are both nontrivial bounds on $\gamma_{4,t}(K), \gamma_{4,s}(K)$, knot invariants which are themselves often difficult to compute. The upper bound is best possible as for instance in the case of $K=4_1$, one has $\gamma_{4,t}(K) = 2=\gamma_{4,s}(K)$ and $u(K)=1$. Also, if $K=8_{18}$ then $\gamma_{4,s}(K) = 3$ while $u(K) = 2$. 

Let $g_{4,t}(K), g_{4,s}(K)$ denote the topological and smooth 4-genus of $K$. It is easy to see that $\gamma_{4,t}(K)\le 2g_{4,t}(K)+1$ and $\gamma_{4,s}(K)\le 2g_{4,s}(K)+1$, while it is well known that $g_{4,s}(K)\le u(K)$. The combination of the two yields $\gamma_{4,t}(K), \gamma_{4,s}(K) \le 2u(K)+1$, substantially weaker bounds than those of Corollary \ref{Corollary on the upper bound of gamma4 by the unknotting number plus one}.   

The lower bound from Corollary \ref{Corollary on the upper bound of gamma4 by the unknotting number plus one} is also best possible as the case of $K=3_1$ shows, for which $\mu_{an}(K) = 1=\gamma_{4,t}(K)=\gamma_{4,s}(K)$. The inequality $\mu_a(K,2)\le u(K)+1$ follows from existing results, see inequality \eqref{Equation 2 with bounds on u_nb} in Section \ref{Section on background on band-unknotting numbers}.

For any slice knot $L$, $\mu_{an}(L)=0$, and in particular, $\mu_{an}\left(K\#(-\bar K)\right)=0$ for any knot $K$. 
\subsection{Organization} The remainder of the paper is organized as follows: Sections \ref{Section with background on various topics} provides background on (dual) band moves, linking forms and band-unknotting numbers. Section \ref{Section with the proof of the first theorem} proves Theorem \ref{Main double inequality for the band-unknotting number} and Corollary \ref{Corollary on the upper bound of gamma4 by the unknotting number plus one}, Section \ref{Section proving strict inequality with respect to max band unknotting number} proves Theorem \ref{Theorem about strict inequality for max band unknotting number}, and Section \ref{Section that proves the theorem about the knot group} provides a proof for Theorem \ref{Theorem about how the knot group does not determine the band unknotting number}. 
\vskip1mm
\noindent {\bf Acknowledgements. } We thank Maggie Miller for a helpful email correspondence.   

\section{Background} \label{Section with background on various topics}
\subsection{Band Moves} \label{Section with background on band moves}
We defined a non-orientable band move on a knot diagram via Figure \ref{Unknotting operations}(b). We elaborate here on this definition, as we shall rely on an equivalent reformulation of it below. 

By a {\em band} we shall mean the oriented boundary of $[0,1]^2$. A {\em band-move} on a knot diagram is a band-surgery, performed by attaching the edges $[0,1]\times \{0,1\}$ of a band to two disjoint arcs of the knot diagram, and then deleting the attaching arcs and adding the edges $\{0,1\}\times [0,1]$ to the diagram. We say that a band-move is an {\em oriented band-move} if the attaching regions of the band have orientation that either both agree or both disagree with the orientations of the attaching arcs, inherited from an orientation of the knot. The knot need not be oriented for this to make sense. One imposes an orientation on the knot just temporarily to check compatibility of orientations when attaching the band, and either orientation on the knot will yield the same compatibility (or non-compatibility) results. 

A {\em non-orientable band move} is one that isn't orientable; it is distinguished by having one attaching edge have the orientation that agrees with the orientation of its attaching arc, with the other edge having the opposite orientation of the attaching arc. This definition of a non-orientable band move agrees with that of Figure \ref{Unknotting operations}(b) as we can always modify a given diagram by planar isotopy and Reidemeister moves, to bring the two attaching arcs of the band close to each other so they look as in Figure \ref{Unknotting operations}(b), and only then perform the band move. The same can of course be done with an oriented band move. We shall take advantage of this broader viewpoint of band moves, and use bands that can meander around and across a given knot diagram between the attaching arcs. 

Some authors refer to orientable band moves as {\em fusion} or {\em fission} moves according to whether they fuse or split components in a knot/link. Similarly, non-orientable band moves are sometimes referred to as {\em half-twists} or {\em crosscaps}. For later use, we make these easy but useful observations. 
\begin{remark} \label{Change of components under band moves}
Let $K$ be an oriented knot and $K_1, K_2$ be two distinct components in an oriented link.  
\begin{itemize} 
\item[(i)] Attaching a non-orientable band to $K$ produces a new knot, while attaching its ends to $K_1$ and $K_2$ fuses the two components into a single component. The new knot/link obtained from a non-orientable band move does not admit an orientation compatible with the original knot/link. 
\item[(ii)] Attaching an orientable band to $K$ produces a two-component link, while attaching its ends to $K_1$ and $K_2$ fuses the two components into a single component. In either case, the resulting knot/link inherits a well-defined orientation from the knot/link it was obtained from. 
\end{itemize} 
\end{remark}

When attaching multiple bands to a diagram, we shall always arrange the attaching arcs of the bands to all be mutually disjoint, and to only lie on the diagram of the knot/link, rather than having attaching arcs of later bands, lie on previously attached bands. 

If $A$ is a band that is attached to the knot/link $K$, we simply write $K\cup A$ for the knot/link obtained by the surgery on $A$. Our notation does not distinguish between orientable and non-orientable band moves. If we attach two bands $A$ and $B$ to $K$, then writing $K\cup A\cup B$ means that we first attach $A$ to $K$ and then $B$ to $K\cup A$, while $K\cup B\cup A$ means that we first attach $B$ to $K$ and then $A$ to $K\cup B$. The final outcomes $K\cup A\cup B$ and $K\cup B\cup A$ are the same (given our convention from the previous paragraph), but the order of attachment may change what kind of band each of $A$ or $B$ is. 
\subsection{Dual Band Moves} \label{Subsection on Dual Band Moves}
Let $K$ be an oriented knot or link and let $A$ be a band, orientable or not. Thus, $A = \partial [0,1]^2$ as an oriented 1-manifold with the edges $[0,1]\times \{0,1\}$ attached to two disjoint arcs of $K$. The orientations of $K$ and of $A$ determine if $A$ is an oriented or non-orientable band attachment. Regardless of its orientability, we call the $\bar A = \partial([0,1]^2)$ the {\em dual band of $A$} if $\bar A$ gets attached to $K\cup A$ via the arcs $\{0,1\}\times [0,1]$. It is easy to see that $K\cup A\cup \bar A = K$ and that the dual band of the dual band, equals the original band. Sometimes, for easier visibility, we shall identify the dual band with $\bar A = \left([0,1]\times [\frac{1}{2}, \frac{3}{4}]\right)$  attached it to $K\cup A$ via $\{0,1\}\times [\frac{1}{2}, \frac{3}{4}]$. In this case $K\cup A\cup \bar A$ is isotopic to $K$ and the dual of the dual is isotopic to the original band. See Figure \ref{Figure on dual band moves} for an illustration. 
\begin{figure}[ht]
\includegraphics[width=0.80\textwidth]{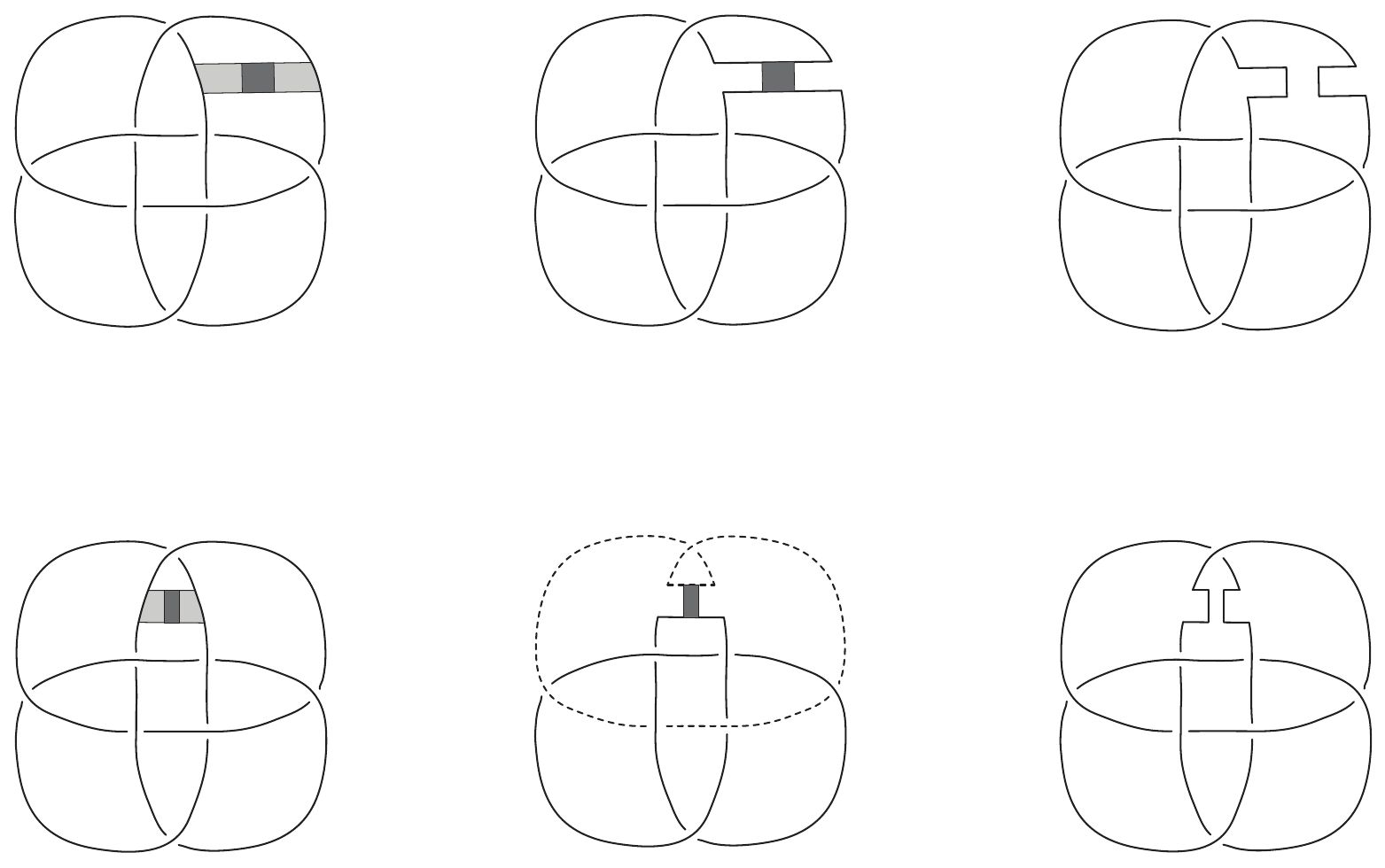}
\put(-315,115){(a)}
\put(-183, 115){(b)}
\put(-50, 115){(c)}
\put(-315,-20){(d)}
\put(-183, -20){(e)}
\put(-50, -20){(f)}
\caption{An illustration of dual band moves. The top row shows a non-orientable band move being performed on the knot $K=8_{18}$. The band itself is indicated in Figure (a) and is shaded in light gray. Its dual band is represented by the darker gray rectangle. In Figure (b) the non-orientable band move has been performed, the dual band remains unchanged. In Figure (c) we perform the dual band move (also a non-orientable band move) to return to an isotopic copy of the original knot from Figure (a). Figures (d)--(f) do the same for an orientable band move. The dotted curve in Figure (e) is to emphasize that a two-component link has been obtained. } \label{Figure on dual band moves}
\end{figure}
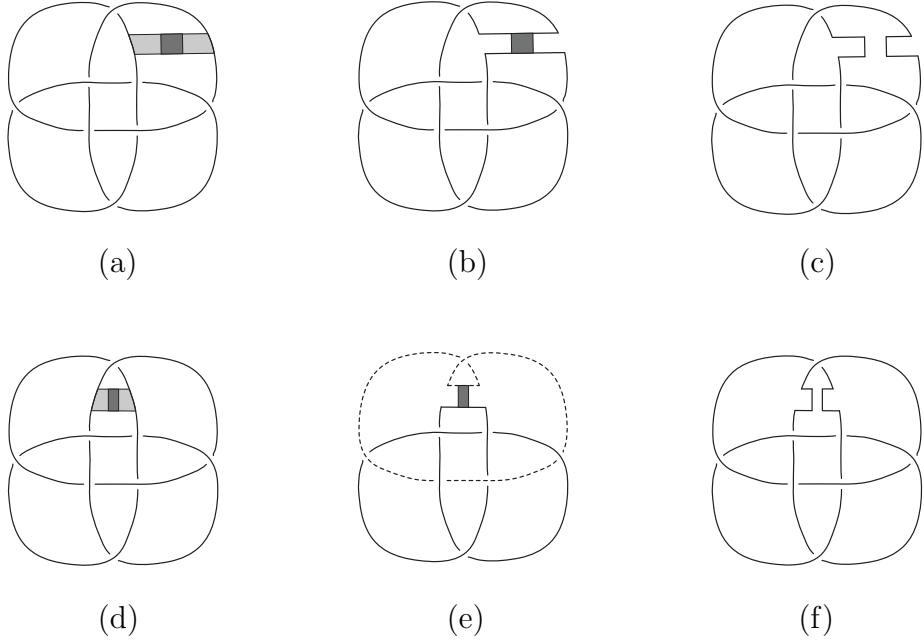
\begin{remark}
Duality of bands preserves the orietability type of the band: The dual band is orientable if and only if the band is orientable; and similarly for non-orientable bands. 
\end{remark}  
\subsection{Linking forms} \label{Section on background on linking forms}
The material presented in this section can be found in most books on bilinear or quadratic forms. We lean on the exposition from \cite[Section 5]{Scharlau}, see also \cite{Kawauchi, Wall}. 

A {\em linking form} is a pair $(G,\lambda)$ consisting of a finite Abelian group $G$ and a symmetric bilinear function $\lambda :G\times G\to \mathbb Q/\mathbb Z$. We sometimes refer to $\lambda$ itself as the linking form if $G$ is understood from context. A linking form $\lambda$ is called {\em non-singular} or {\em regular} if its dual 
$$\lambda^\ast :G\to G^\ast = \textrm{Hom}(G, \mathbb Q/\mathbb Z), \quad \text{ define as } \quad \lambda^\ast (g) = \lambda (g, \cdot),$$
is an isomorphism. Two linking forms $(G_1, \lambda_1)$ and $(G_2, \lambda_2)$ are called {\em isomorphic} or {\em isometric} if there exists an isomorphism $\varphi:G_1\to G_2$ of groups such that $\lambda_2(\varphi(g_1), \varphi(g_2)) = \lambda _1(g_1, g_2)$ for all $g_1, g_2\in G_1$. This is an equivalence relation among linking forms, and we shall without further mention identify forms that are isomorphic. 

Two linking forms $(G_1, \lambda_1)$ and $(G_2, \lambda_2)$ can be added, and we write $(G_1, \lambda_1)\oplus (G_2, \lambda_2)$ for this sum, by considering the function $\lambda_1\oplus \lambda_2$ on $(G_1\oplus G_2)\times (G_1\oplus G_2)$ defined by $\lambda_1\oplus \lambda_2((g_1,h_1),(g_2,h_2)) = \lambda_1(g_1,g_2) + \lambda_2(h_1,h_2)$. The sum of non-singular linking forms is again non-singular. 

For a linking form $(G,\lambda)$ and a subgroup $H$ of $G$, we define its {\em orthogonal complement} $H^\perp$ as $H^\perp = \{g\in G\,|\, \lambda (g,h)=0, \text{ for all } h\in H\}$. A key result about non-singular linking forms (cf. \cite[Lemma 1]{Wall}) is the orthogonal splitting 
\begin{equation} \label{Splitting of linking form by orthogonal complement}
(G, \lambda) \cong (H,\lambda_1)\oplus (H^\perp, \lambda_2),
\end{equation}
valid under the hypothesis that $\lambda_1: = \lambda|_{H\times H}$ is   non-singular. In this case, $\lambda _2:= \lambda|_{H^\perp\times H^\perp}$ is also non-singular. 

A form $(G,\lambda)$ is called {\em isotropic} if there exists an {\em isotropic element $g\in G$}, that is a non-zero element with $\lambda(g,g)=0$; otherwise it is called {\em anisotropic}. If $\lambda$ is non-singular and isotropic and $g$ is an isotropic element, one can always find an element $h\in G$ with $\lambda(g,h)\ne0$ and such that the restriction of $\lambda$ to $H\times H$ is non-singular. Here $H$ is the subgroup of $G$ generated by $g$ and $h$. By \eqref{Splitting of linking form by orthogonal complement}, this leads to the isomorphism 
\begin{equation} \label{Equation with decomposition of non-singular linking forms}
(G,\lambda) \cong (H, \lambda|_{H\times H})\oplus (H^\perp,\lambda|_{H^\perp\times H^\perp}).	
\end{equation}
Clearly $H^\perp$ is a proper subgroup of $G$, and if $H^\perp$ is also isotropic, on can continue this process of splitting off summands generated by isotropic elements, until one arrives either at an anisotropic linking form or the trivial subgroup of $G$ (which is by convention also considered anisotropic). This sketches the proof of the next theorem, known as the Witt Decomposition Theorem \cite[5.11 Corollary]{Scharlau}. To state it, recall that a linking form $(H,\mu)$ is called {\em metabolic} if $H$ has a subgroup $K$ with $K=K^\perp$. The form $(H,\lambda|_{H\times H})$ in decomposition \eqref{Equation with decomposition of non-singular linking forms} is metabolic. 
\begin{utheorem}[Witt Decomposition]
Let $(G,\lambda)$ be a non-singular linking form. Then there exists an isomorphism 
$$(G,\lambda) \cong (G_1, \lambda_1)\oplus (G_2, \lambda_2),$$
with $(G_1, \lambda_1)$ anisotropic (possibly zero) and with $(G_2, \lambda_2)$ metabolic. The isomorphism type of $(G_1, \lambda_1)$ is uniquely determined by $(G, \lambda)$. 
\end{utheorem} 
The form $(G_1, \lambda_1)$ from the previous theorem is called the {\em anisotropic component} or {\em anisotropic part} or the {\em kernel form} of $(G,\lambda)$. 
\vskip2mm
For a knot $K$, let $Y = \Sigma_2(K)$ by the cyclic two-fold cover of $S^3$ branched along $K$. Then $Y$ is a rational homology $3$-sphere with $|H_1(Y)| =\det K$. The linking form $\lambda :H_1(Y)\times H_1(Y)\to \mathbb Q/\mathbb Z$ is defined as follows: Let $a,b\in H_1(Y)$ be two homology classes represented by embedded curves $\alpha, \beta\subset Y$. Let $m\in \mathbb N$ be such that $m\cdot b = 0 \in H_1(Y)$, and thus $m\cdot \beta =\partial \sigma$ for some 2-chain $\sigma\subset Y$. Then 
$$\lambda(a,b) = \frac{1}{m}\, \alpha \cdot \sigma,$$
where on the right-hand, $\alpha\cdot \sigma$ is the linking pairing that counts intersection points between $\alpha$ and $\sigma$ with sign. This linking pairing is well known to be non-singular.    
\vskip1mm
Let $X$ be a compact, connected, oriented 4-manifold $X$ with $\partial X = Y$. 
The long exact homology sequence for the pair $(X,Y)$ includes these terms:
\begin{align} \label{Homology exact sequence}
\dots H_2(X) \stackrel{\iota_\ast}{\longrightarrow} H_2(X,Y)\stackrel{\partial}{\longrightarrow} H_1(Y) \to H_1(X) \to \dots  
\end{align}
Write $H_2(X)\cong F\oplus T_1$ and $H_2(X,Y) \cong F\oplus T_2$ with $T_1, T_2$ the torsion subgroubs of $H_2(X), H_2(X,Y)$ and with $F$ a free Abelian group of rank $b_2(X)$. The component of the map $H_2(X)\to H_2(X,Y)$ in the above exact sequence that maps $F\subset H_2(X)$ to $F\subset H_2(X,Y)$, can, with a suitable choice of bases, be represented by a non-singular, symmetric, integral matrix $A$, of dimension $b_2(X)\times b_2(X)$. 

Let $B=\partial (F)\subset H_1(Y)$. By the exact sequence \eqref{Homology exact sequence}, $B\cong \cokerr(A)$, and the minimum number of generators of the cokernel of $A$, is given by the number of diagonal entries in the Smith Normal Form of $A$ that are strictly greater than $1$. In particular, the minimum number of such generators is no greater than $b_2(X)$, showing that the minimum number of generators of $B$ is also bounded from above by $b_2(X)$.   

It was shown in \cite{Gilmer} that $\lambda|_{B\times B}$ is non-singular, and so by \eqref{Splitting of linking form by orthogonal complement} one can decompose the linking form $\lambda$  as 
$$(H_1(Y), \lambda) \cong (B, \lambda _0)\oplus (B^\perp, \lambda_3), $$
with $\lambda _0 = \lambda|_{B\times B}$ and $\lambda_3 = \lambda_{B^\perp \times B^\perp}$. \cite[Lemma 1]{Gilmer} proves that $(B^\perp,\lambda_3)$ is metabolic. Apply the Witt Decomposition to $(B, \lambda_0)$ to obtain 
$$(B,\lambda_0) \cong (G_1, \lambda_1)\oplus (G_2, \lambda_2)$$
with $(G_1,\lambda_1)$ anisotropic and with $(G_2,\lambda_2)$ metabolic. It is easy to see that the orthogonal direct sum of two metabolic linking forms is again metabolic, yielding this Witt Decomposition of $(H_1(Y),\lambda)$:
\begin{equation} \label{With decomposition of the linking form on the two-fold cyclic branched cover}
(H_1(Y), \lambda) \cong (G_1, \lambda_1) \oplus \left[(G_2, \lambda_2) \oplus (B^\perp, \lambda_3) \right].
\end{equation}
By 
Recall that $\mu_{an}(K)$ is the minimum number of generators of $G_1$. Since $G_1\le G_1\oplus G_2 \cong B \cong \cokerr(A)$, we have proved this statement.  
\begin{theorem} \label{Dragon Kung Fu Theorem}
Let $K$ be a knot and let $X$ be a compact, connected, and oriented 4-manifold with $\partial X = \Sigma_2(X)$. Then  
$$
\mu_{an}(K)  \le b_2(X). 
$$
\end{theorem}

\subsection{The band-unknotting number} \label{Section on background on band-unknotting numbers}
\begin{figure}[ht]
	\includegraphics[width=0.70\textwidth]{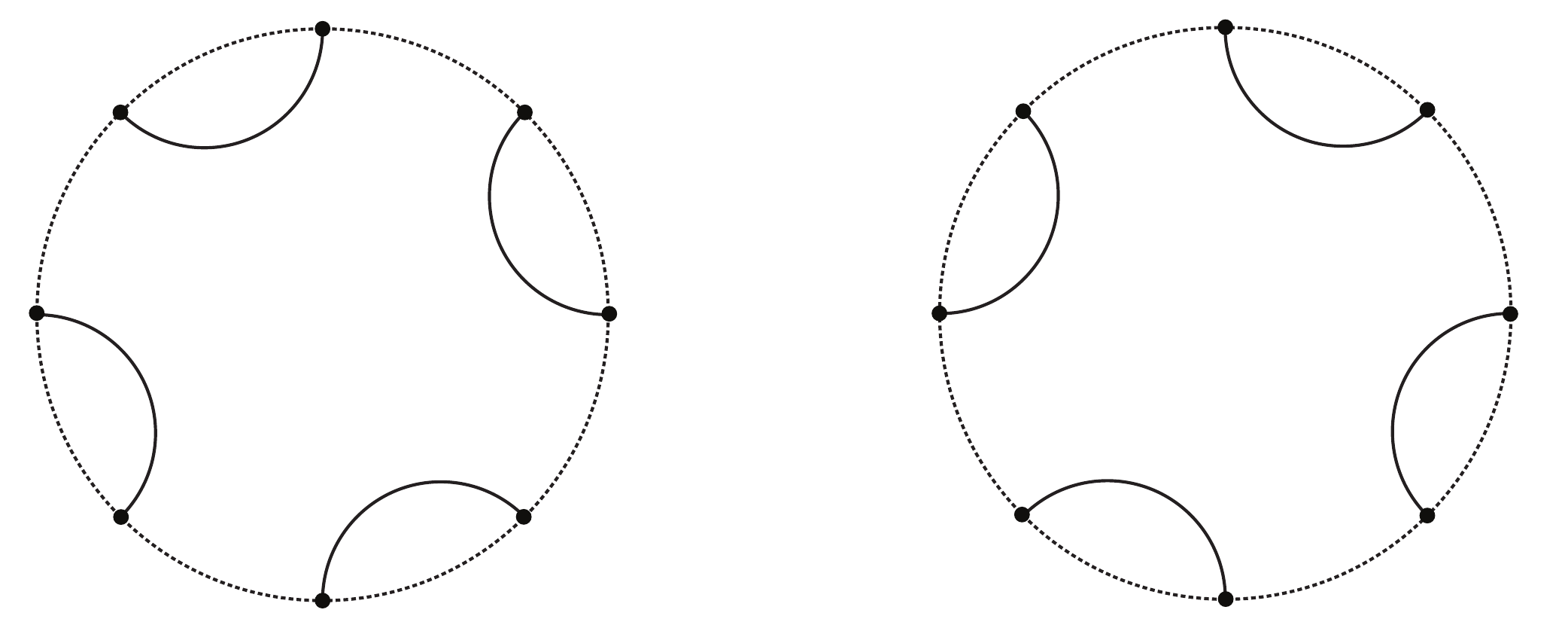}
	\put(-253,-15){(a)}
	\put(-75, -15){(b)}
	\caption{(a) The unknotting operation $H(n)$, $n\ge 2$, shown here for $n=4$. Under this move, a knot/link diagram is locally modified by replacing the arrangement of $n$ unknotted arcs as in Figure (a), with the arrangement as in Figure (b). One requires that the number of components be preserved. An $H(2)$-move is a non-orientable band move. } \label{H(n) move}
\end{figure}

The band-unknotting number as we have defined it in the introduction, has been considered by other authors, sometimes under different names. The earliest result about band-unknotting that we are aware of goes back to Lickorish \cite{Lickorish} from 1986. What we call a non-orientable band move, Lickorish refers to ``adding a twisted band'' to a knot diagram, and he proves that the Figure Eight knot cannot be rendered unknotted by adding a single twisted band, a claim that follows from this main result of \cite{Lickorish}:
\begin{utheorem}[Lickorish \cite{Lickorish}] If $K$ is a knot with $u_{nb}(K)=1$, then $H_1(\Sigma_2(K))$ is cyclic and there exists a generator $\alpha\in H_1(\Sigma_2(K))$ such that  
\begin{equation} \label{Equation 1 with bounds on u_nb}
\lambda (\alpha, \alpha) = \pm \frac{1}{\det K}.
\end{equation} 
\end{utheorem}
\noindent It follows that if $K$ is any knot with $\det K>1$, then $u_{nb}(K\#(\pm K)), u_{nb}(K\# \pm \bar K) \ge 2$. While Lickorish's work is related to band-unknotting, he does not explicitly mention it. 
\vskip1mm
Hoste, Nakanishi and Taniyama \cite{HosteNakanishiTaniyama} introduced a family of unknotting operations $H(n)$, indexed by integers $n\ge 2$, described in Figure \ref{H(n) move}.  
The corresponding unknotting numbers are denoted by $u_n(K)$, and it is easy to check that an $H(2)$-move equals a non-orientable band move, so that $u_{nb}(K) = u_2(K)$. To quote results from their paper we shall rely on below, let $\Sigma_r(K)$ denote the $r$-fold cyclic cover of $S^3$ branched along $K$, with $r\ge 2$. Let $\mu(K,r)$ denote the minimum number of generators of $H_1(\Sigma_r(K))$, and for simplicity, let $\mu(K) = \mu(K,2)$. 
\begin{utheorem} \cite[Theorem 4]{HosteNakanishiTaniyama} For all integers $r\ge 2$, one obtains 
\begin{equation} \label{Equation 2 with bounds on u_nb}
\frac{\mu(K,r)}{r-1} \le u_{nb}(K),\quad \text{ and thus } \quad \mu(K) \le u_{nb}(K). 
\end{equation} 
\end{utheorem} 
\vskip1mm
Kanenobu and Miyazawa \cite{KanenobuMiyazawa} study $u_2(K)=u_{nb}(K)$ further, proving various boundsand computing its value for most knots with 9 or fewer crossings. We shall make use of this result:
\begin{utheorem} \cite[Theorem 3.1]{KanenobuMiyazawa} 
 The Gordian and band-unknotting numbers are related by the inequality:
\begin{equation} \label{Equation 4 with bounds on u_nb}
u_{nb}(K) \le u(K)+1.
\end{equation}

\end{utheorem}
\vskip1mm
Yasuhara \cite{Yasuhara} showed that if $K$ is a knot with $\sigma(K) + 4\cdot \mathrm{Arf}(K)\equiv 0\pmod 8$, then $\gamma_{4,s}(K)\ge 2$. It follows that such knots must satisfy the inequality $u_{nb}(K)\ge 2$, though we shall not make use of this here.    
\section{The proof of Theorem \ref{Main double inequality for the band-unknotting number} and Corollary \ref{Corollary on the upper bound of gamma4 by the unknotting number plus one}} \label{Section with the proof of the first theorem}
\subsection{Part (a) of Theorem \ref{Main double inequality for the band-unknotting number}} Part (a) of Theorem \ref{Main double inequality for the band-unknotting number} claims the double inequality 
$$ \mu_{an}(K) \le \gamma_4(K) \le u_{nb}(K),$$
with $\gamma_4(K)$ equal to either $\gamma_{4,t}(K)$ or $\gamma_{4,s}(K)$. 
The bound $\gamma_4(K) \le u_{nb}(K)$ follows from two facts established in \cite{JabukaKelly}:
\begin{itemize}
\item[(i)] If $K$ and $K'$ are two knots that possess diagrams that are related by a single non-orientable band move, then $|\gamma_4(K) - \gamma_4(K')|\le 1$.  
\item[(ii)] If $K$ is a non-trivial knot related to a slice knot by a single non-orientable band move, then $\gamma_4(K) = 1$. 
\end{itemize}
The work in \cite{JabukaKelly} focuses on $\gamma_{4,s}$, but the two statements above are readily verified to be valid for $\gamma_{4,t}$ as well. 

If $K$ is a knot with $\gamma_4(K) = k\ge 3$, then after $k-2$ non-orientable band moves on a diagram representing $K$, we obtain a diagram representing a knot $K'$ with $\gamma_4(K')\ge 2$. If an additional band move on $K'$ led to the unknot, that would erroneously imply $\gamma_4(K')=1$, showing that at least two more non-orientable band moves are needed to unknot $K'$, and thus a total of at least $k=\gamma_4(K)$ band moves to unknot $K$. This proves Part (a) of Theorem \ref{Main double inequality for the band-unknotting number}. 
\vskip1mm
For the lower bound $\mu_{an}(K)\le \gamma_4(K)$ from Theorem \ref{Main double inequality for the band-unknotting number}, let $\sigma\subset D^4$ be a properly embedded, non-orientable surface with $\partial \sigma =K$. Assume that either $\sigma$ is locally topologically flatly embedded and that $\gamma_{4,t}(K)=b_1(\sigma)$, or else that $\sigma$ is smoothly embedded with $\gamma_{4,s}(K)=b_1(\sigma)$. In either case, let $X$ be the connected, compact, oriented 4-manifold obtained as the 2-fold cyclic cover of $D^4$ with branching set $\sigma$. Then $b_2(X) = b_1(\sigma)$ and $X$ meets the hypotheses of Theorem \ref{Dragon Kung Fu Theorem}, showing that $\mu_{an}(K)\le b_2(X)=\gamma_4(K)$. This completes the proof of Part (a) of Theorem \ref{Main double inequality for the band-unknotting number}.
\vskip3mm
\subsection{Part (b) of Theorem \ref{Main double inequality for the band-unknotting number}} Endow $K$ with a choice of orientation. Let $n=F\left(K\#(-\bar K)\right)$ so that by definition of the fusion number, the knot $K\#(-\bar K)$ is the boundary of a ribbon disk $D^2\looparrowright S^3$ that becomes the $(n+1)$-component oriented unlink $U_{n+1}$ (the orientation induced by that on $K$) after applying $n$ orientable band moves $A_1, \dots, A_n$ to $K\#(-\bar K)$: 
$$U_{n+1} = \left(K\#(-\bar K)\right) \cup A_1\cup \dots \cup A_n.$$
By Remark \ref{Change of components under band moves}, we can only obtain an $(n+1)$-component link $U_{n+1}$ from a knot by using $n$ oriented band moves, if each band has both ends attached to the same component. Accordingly, each oriented band move $A_i$ splits its attaching component into two components $C_i$ and $C'_i$, and the dual band $\bar A_i$ has its two attaching arcs lying on $C_i$ and $C'_i$, thereby fusing two components with an orientable band move. 

We may isotope $U_{n+1}$ to the standard diagram of $U_{n+1}$, exercising care to ``drag'' the dual bands $\bar A_1, \dots, \bar A_n$ along with the isotopy. In this manner, we obtain a configuration of the standard unlink connected by the oriented bands $\bar A_1, \dots, \bar A_n$, and attaching the said dual bands, returns us to the original knot $K\#(-\bar K)$. 
\begin{proposition} \label{Proposition with all the dual band moves}
Let $K$ be an oriented knot with $n=F\left(K\#(-\bar K)\right)$ and suppose that the configuration of the unlink $U_{n+1}$ and the dual bands $\bar A_1, \dots, \bar A_n$ was obtained from from $K\#(-\bar K)$ by attaching the $n$ orientable bands $A_1, \dots, A_n$. Then there exist $n$ non-orientable band moves $B_1, \dots, B_n$ on $U_{n+1}$ such that:
\begin{itemize}
\item[(i)] $U_{n+1}\cup B_1\cup \dots \cup B_n = U_1=$unknot. 
\item[(ii)] In the diagram for $U_1$ obtained in Part (i), the band moves $\bar A_1, \dots, \bar A_n$ are non-orientable band moves. 
\end{itemize}
\end{proposition}
The preceding proposition implies that  
\begin{align*}
	U_{n+1}  & =  \left( K\#(-\bar K)\right) \cup A_1 \cup \dots \cup A_n \cr 
	U_{n+1}\cup B_1\cup \dots \cup B_n  & =  \left( K\#(-\bar K)\right)\cup A_1 \cup \dots \cup A_n\cup B_1\cup \dots \cup B_n \cr 
	U_1& =  \left(K\#(-\bar K)\right)\cup B_1\cup \dots \cup B_n \cup A_1 \cup \dots \cup A_n,
\end{align*}
with all the $2n$ band moves on the right-hand side of the last line being non-orientable band moves. Thus $u_{nb}(K\#(-\bar K))\le 2n$, proving Theorem \ref{Main double inequality for the band-unknotting number}. It remains to proof the proposition.
\begin{proof} (Of Proposition \ref{Proposition with all the dual band moves}). We proceed by induction on $n\ge 1$.  
\vskip1mm
\noindent {\bf Base of induction.} The simple case of $F\left( K\#(-\bar K)\right)=1$ is worked out in Figure \ref{Band_moves_trading_induction_base}. 
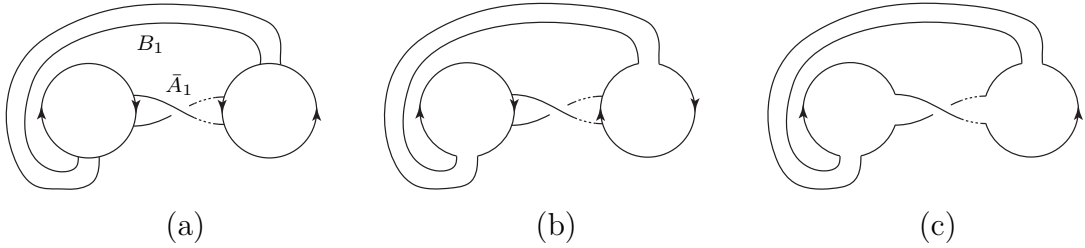
\begin{figure}[ht]
\includegraphics[width=0.95\textwidth]{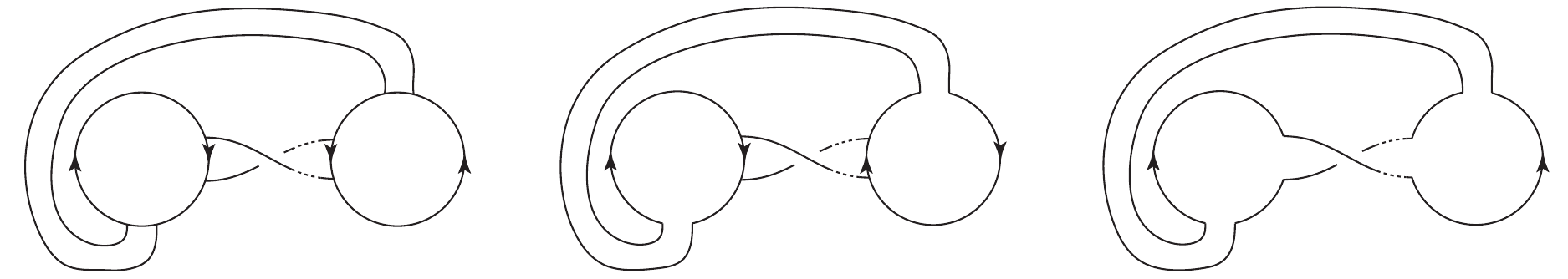}
\put(-352,-15){(a)}
\put(-352,39){\tiny $\bar A_1$}
\put(-363,55){\tiny $B_1$}
\put(-212,-15){(b)}
\put(-67, -15){(c)}
\caption{Figure (a) shows an oriented 2-component unlink $U_2$ isotoped into standard position. The single dual band $\bar A_1$ is the half-twisted band connecting the two components of $U_2$, and were we to first perform this band move, it would be an oriented band move. The dotted lines in this band indicate that the band $\bar A_1$ can be complicated; it can be knotted and link the components of $U_2$. The band $B_1$ is the untwisted band on the outside of both components of $U_2$. In Figure (b) we performed the non-orientable band move $B_1$, and adjusted the orientation on the thus obtained unknot $U_1$. In Figure (c) we finally perform the band move $\bar A_1$ which is now a non-orientable band move.  } \label{Band_moves_trading_induction_base}
\end{figure}
The conclusion from Figure \ref{Band_moves_trading_induction_base} is that $B_1$ is a non-orientable band move on $U_2$ and that $U_2\cup B_1$ the unknot $U_1$. The band move $\bar A_1$ performed on $U_2\cup B_1$, is a non-orientable band move. We have thus shown that 
\begin{align*}
U_2 &= \left(K\#(-\bar K)\right)\cup A_1 \cr
U_2 \cup B_1 & = \left(K\#(-\bar K)\right)\cup A_1 \cup B_1, \cr 
U_1 & = \left(K\#(-\bar K)\right)\cup B_1 \cup A_1.  
\end{align*}
\vskip1mm 
\noindent {\bf Step of induction.} Assume the proposition has been verified for fusion numbers no greater than $n-1$, and let $K$ be an oriented knot with $F\left(K\#(-\bar K)\right)=n$. Perform $n$ oriented band moves $A_1, \dots, A_n$ on $K\#(-\bar K)$ to obtain the oriented $(n+1)$-component unlink $U_{n+1}$. Pick a component $C$ of $U_{n+1}$ that is only connected to one single other component of $U_{n+1}$ by a single dual band, say $\bar A_n$ for concreteness. This is always possible for if not, then every component would be connected to at least two other components, and would thus have at least two dual bands attached to it, making the total number of dual bands exceed $\frac{1}{2}\cdot (2(n+1)) = n+1$. But there are only $n$ bands, a contradiction. 

Thus we can always arrange to be in the situation of Figure \ref{Band_moves_trading_induction_step}, which we explain now. Since none of the bands $A_1, \dots, A_{n-1}$ attach to the component $C$, then $\left((K\#(-\bar K)\right) \cup A_n = L\cup C$, with $L$ some knot with orientation induced by that of $K\#(-\bar K)$. Attaching the oriented bands $A_1, \dots, A_{n-1}$ to $L$ yields the unlink $U_n$ and thus, by induction hypotheses, there exist non-orientable bands $B_1, \dots, B_{n-1}$ satisfying the two conclusions from the statement of Proposition \ref{Proposition with all the dual band moves}. The result of performing the non-orientable band moves $L\cup B_1\cup \dots B_{n-1}\cup \bar A_1\cup \dots \cup \bar A_{n-1}$ yields a knot $L' = L\cup B_1\cup \dots \cup B_{n-1}$, oriented in an arbitrary fashion. The knot $L'$ is shown in the ``black box'' in each of Parts (a)--(c) in Figure \ref{Band_moves_trading_induction_step}. Only an unknotted arc protrudes past the box, used to attach the dual band $\bar A_n$ between $L'$ and $C$. 
\begin{figure}[ht]
\includegraphics[width=0.60\textwidth]{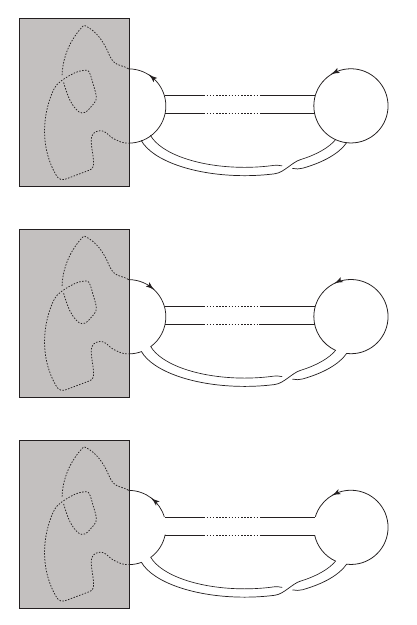}
\put(-290,330){(a)}
\put(-123,340){$\bar A_n$}
\put(-123,296){$B_n$}
\put(-250,283){$L'$}
\put(-25,303){$C$}
\put(-290,195){(b)}
\put(-123,207){$\bar A_n$}
\put(-123,162){$B_n$}
\put(-250,150){$L'$}
\put(-25,170){$C$}
\put(-290,60){(c)}
\put(-123,72){$\bar A_n$}
\put(-123,28){$B_n$}
\put(-250,17){$L'$}
\put(-25,37){$C$}
\caption{The inductive step in the proof of Proposition \ref{Proposition with all the dual band moves}. Figure (a) shows the knot $L' = U_n\cup B_1\cup \dots \cup B_{n-1}\cup \bar A_1\cup \dots \cup \bar A_{n-1}$ and the unknotted component $C$. Indicated are also the bands $\bar A_n$ and $B_n$, though no band moves on either have been performed. Figure (b) shows the result of performing the non-orientable band move $B_n$, while Figure (c) shows the outcome of performing the non-orientable band move $\bar A_n$.  } \label{Band_moves_trading_induction_step}
\end{figure}
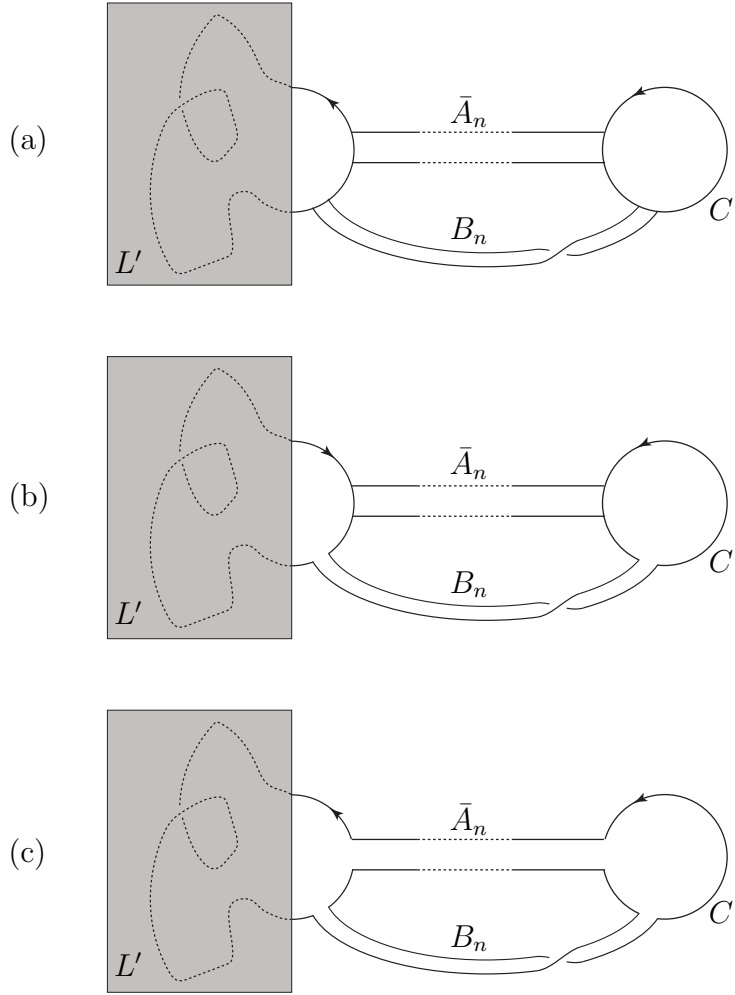
The dots on the band $\bar A_n$ indicate again that this band may be very complicated: It can itself be knotted, it can be linked with $C$, and even enter and exit the black box multiple times. The non-oriented band $B_n$ is also shown. In Part (a) of Figure \ref{Band_moves_trading_induction_step}, no band moves on either $\bar A_n$ or $B_n$ have been performed, but observe that if we were to perform the move on $\bar A_n$ first, we would be performing an oriented band move. In Part (b) we performed the non-orientable band move on band $B_n$, while finally in Part (c), we perform what is now a non-orientable band move on $\bar A_n$. The final result, as seen in Part (c), is 
$$L'\cup B_n \cup \bar A_n = U_n\cup B_1\cup \dots \cup B_{n-1} \cup \bar A_1\cup \dots \cup A_{n-1}\cup B_n \cup \bar A_n.$$ 
As already explained above (in the paragraph following the statement of Proposition \ref{Proposition with all the dual band moves}), this proves Part (ii) of the proposition. Part (i) is also clear for $U_n\cup B_1\cup\dots \cup B_{n-1}$ is the unknot, by induction hypothesis, and clearly so is $U_n\cup B_1\cup\dots \cup B_{n-1}\cup B_n$. 
\end{proof}
\subsection{Proof of Corollary \ref{Corollary on the upper bound of gamma4 by the unknotting number plus one}}
Corollary \ref{Corollary on the upper bound of gamma4 by the unknotting number plus one} is a direct consequence of Theorem \ref{Main double inequality for the band-unknotting number} and the second inequality in \eqref{Equation 2 with bounds on u_nb}.
\section{Proof of Theorem \ref{Theorem about strict inequality for max band unknotting number}} \label{Section proving strict inequality with respect to max band unknotting number}
\begin{figure}[ht]
\includegraphics[width=0.50\textwidth]{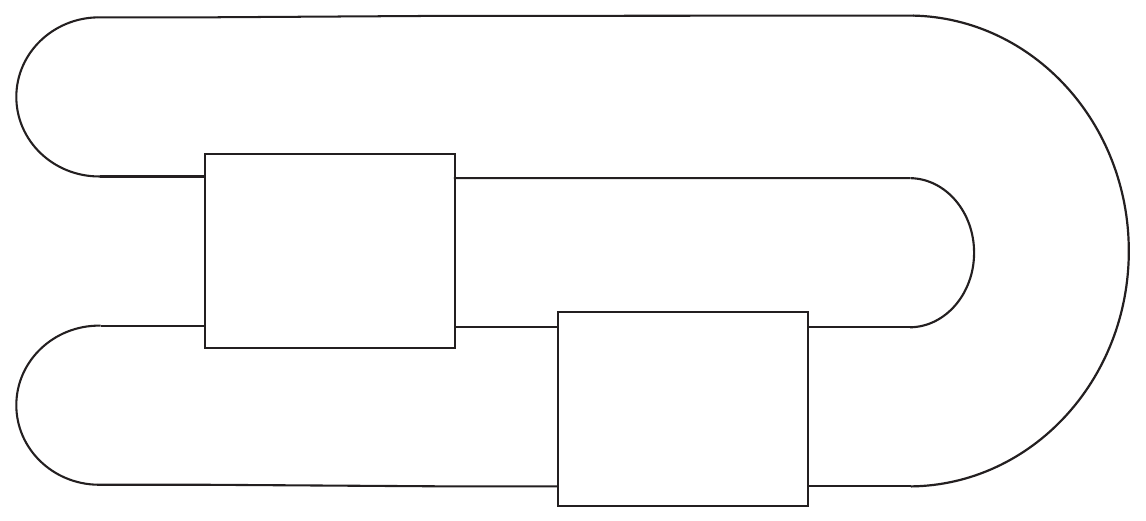}
\put(-163,50){$m$}
\put(-100,18){$-n$}
\caption{The {\em double twist knot} $C(m,n)$ is the two-bridge knot with two twisting regions. Each rectangle labeled by an integer represents a pair of parallel strands with with as many half-twists between them as indicated by the absolute value of the integer. A positive integer indicates right-handed half-twists, while a negative integer indicates left-handed half-twists. The determinant of $C(m,n)$ is $\left|\frac{mn+1}{n}\right|$, showing that $C(m,n)$ is a knot if and only if at least one of $m$ or $n$ is even.} \label{Figure on double twist knots}
\end{figure}

In \cite{HosteShanahanVanCott}, the smooth non-orientble 4-genus $\gamma_{4,s}$ of the family of ``double twists'' $C(m,n)$ (Figure \ref{Figure on double twist knots}) is studied. 
All double twist knots are two-bride knots and thus br$(C(m,n))=2$. On the other hand, Theorem 6.2 in \cite{HosteShanahanVanCott} proves that 
$$\gamma_{4,s}(C(22+8k, 62+8k))=3, \qquad k\ge 0.$$
It follows that if we let $K_1 = C(22+8k, 62+8k)$ and $K_2= -\bar K_1$ for any of the infinitely many choices of $k\ge 0$, then 
\begin{align*}
u_{nb}(K_1\#K_2) & \le 2 F(K_1\#K_2) \le 2(\text{br}(K_1)-1) = 2, \quad \mathrm{and,}\cr 
3 = \gamma_{4,s}(K_i) & \le u_{nb}(K_i), \quad i=1,2.
\end{align*} 
By Lickorish's results \eqref{Equation 1 with bounds on u_nb}, we obtain $u_{nb}(K_1\#K_2) = 2$. It follows that 
$$2= u_{nb}(K_1\#K_2) < 3 = u_{nb}(K_i), i=1,2,$$
proving Theorem \ref{Theorem about strict inequality for max band unknotting number}. We don't know the band-unknotting number of $C(22+8k, 62+8k)$ but it is not hard to show that $u_{nb}(C(22+8k,62+8k))\le 7+2k$. 
\section{Comments on Theorem \ref{First theorem on inequality among band-unknotting numbers}} \label{Section with comments on the proof of the first main theorem}
Theorem \ref{First theorem on inequality among band-unknotting numbers} follows from Theorem \ref{Theorem about strict inequality for max band unknotting number}, which we proved in Section \ref{Section proving strict inequality with respect to max band unknotting number}. We want to point out that it can be verified independently. Namely, given any two-bridge knot $K$ with $\det K>1$ and $\gamma_{4,s}(K)=2$, it follows that $F\left(K \#(-\bar K)\right)=1$ and thus that $u_{nb}(K\#(-\bar K))=2$. On the other hand, $u_{nb}(K)=u_{nb}(-\bar K) \ge \gamma_{4,s}(K)\ge 2$. Any such knot leads to the strict inequality $u_{nb}(K\#(-\bar K)) < u_{nb}(K)+u_{nb}(-\bar K)$. There are infinitely examles of such two-bridge knots that have been identified in \cite[Theorem 5.8]{HosteShanahanVanCott}.    
\section{Proof of Theorem \ref{Theorem about how the knot group does not determine the band unknotting number}} \label{Section that proves the theorem about the knot group}
\begin{figure}[ht]
\includegraphics[width=0.80\textwidth]{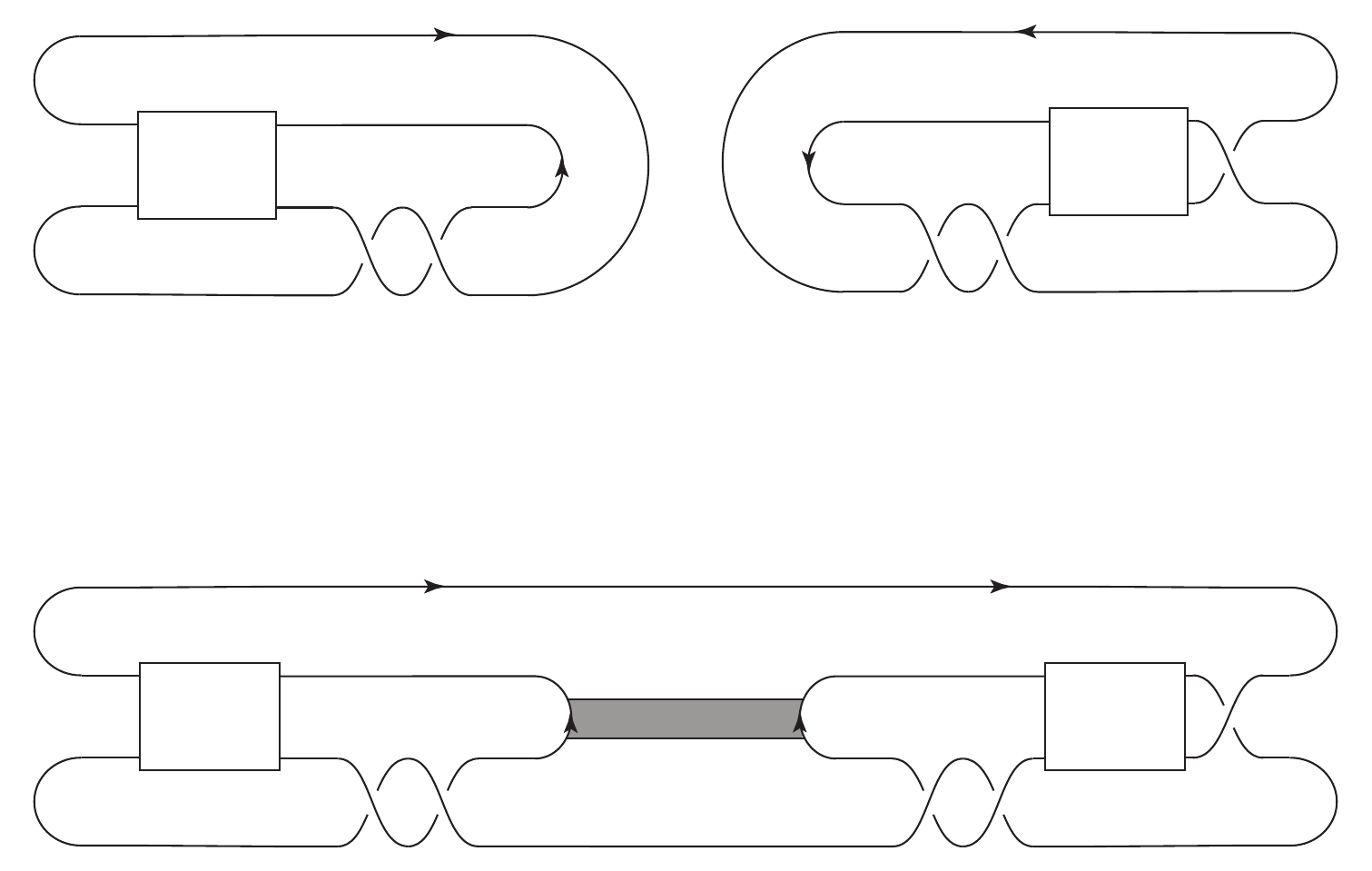}
\put(-302,179){$a$}
\put(-68,179){$a$}
\put(-302,36){$a$}
\put(-75,36){$-a$}
\put(-270,127){(a)}
\put(-100,127){(b)}
\put(-190,-16){(c)}
\caption{For an integer $a\ge 4$, 
Figure (a) shows the two-bridge knot $K_1=K_{(2a-1)/2}$ and Figure (b) the two-bridge knot $K_2=K_{(2a+1)/2}$. Figure (c) show that connected sum $K_1\#(-\bar K_2)$ and indicates the single non-orientable band move that unknots it, showing that $u_{nb}(K_1\#(-\bar K_2)=1$.} \label{Connected sum of two two-bridge knots}
\end{figure}
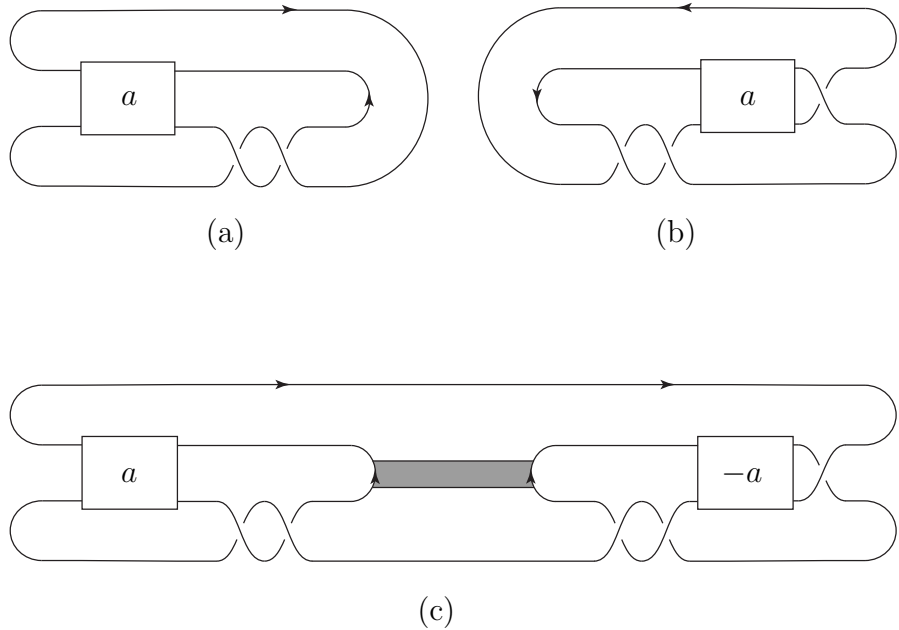

The proof of Theorem \ref{Theorem about how the knot group does not determine the band unknotting number} is inspired by \cite[Theorem 9.1]{KanenobuMiyazawa}. 
For a positive integer $a\ge 4$ and $a\equiv 4 \pmod 6$ consider the two-bridge knots $K_1$ and $K_2$ whose two-fold branched covers are the lens spaces $L(2a-1,2)$ and $L(2a+1,2)$. The linking forms $\lambda_1$ and $\lambda_2$ on $\Sigma_2(K_1)$ and $\Sigma_2(K_2)$ are 
$$\lambda_1=\left\langle \frac{2}{2a-1} \right\rangle \qquad \text{ and } \qquad \lambda_2=\left\langle \frac{2}{2a+1}\right\rangle. $$
Note that $\gcd(2a-1, 2a+1)=1$ and thus the first homology of the two-fold branched cover of $K_1\#K_2$ is given by $H_1(\Sigma_2(K_1\#K_2))\cong \mathbb Z_{2a-1}\oplus \mathbb Z_{2a+1} \cong \mathbb Z_{4a^2-1}$, and its linking form $\lambda _3\cong \lambda_1\oplus \lambda_2$ is determined by 
$$\lambda_3(1,1) = \left\langle \frac{2a}{4a^2-1}\right\rangle.$$
If we had $u_{nb}(K_1\#K_2)=1$ then by \eqref{Equation 1 with bounds on u_nb} we would obtain an isomorphism $\lambda_3 \cong \left\langle \pm \frac{1}{4a^2-1}\right\rangle$. This would imply the existence of an integer $\lambda$ solving the congruence 
$$2a \equiv \pm \lambda^2 \pmod{(4a^2-1)}.$$
Since $4a^2-1=(2a-1)(2a+1)$, and the two factors are relatively prime, then the above congruence is equivalent to the two congruence relations
$$2a \equiv \pm \lambda^2 \pmod{(2a-1)}\quad \text{ and } \quad  2a\equiv \pm \lambda^2 \pmod{(2a+1)}.$$
Furthermore, because $2a = (2a\pm 1) \mp 1 \equiv \mp 1 \pmod{(2a\pm 1)}$, the above can be rewritten more succinctly as
\begin{equation} \label{Equation with the pair of congruence relations}
1 \equiv \pm \lambda^2 \pmod{(2a-1)}\quad \text{ and } \quad  -1\equiv \pm \lambda^2 \pmod{(2a+1)}.
\end{equation}
Recall now that $a=4+6k$ for $k\ge 0$, and thus $2a+1\equiv 0 \pmod 3$. 
\vskip1mm
\noindent {\bf Case (i). } Suppose the sign in front of $\lambda^2$ is a plus sign. Then the mod 3 version of the second equation in  \eqref{Equation with the pair of congruence relations} reads $-1\equiv \lambda^2 \pmod 3$. But $-1$ is not a square in $\mathbb Z_3$ and so no solution $\lambda$ can exist. \\
{\bf Case (ii). } If the sign in front of $\lambda^2$ equals $-1$, then the first equation in \eqref{Equation with the pair of congruence relations} becomes $-1\equiv \lambda^2 \pmod{(2a-1)}$, a congruence saying that $-1$ is a quadratic residue modulo $2a-1$. Determining when $-1$ is a quadratic residue is relatively easy. Recall that for an odd integer $n\ge 3$ and a coprime non-zero integer $x$, the Jacobi symbol $\left(\frac{x}{n}\right)\in \{\pm 1\}$ equals 1 when $x$ is a quadratic residue modulo $n$. The key fact for us, see for example \cite[Theorem 11.10]{Rosen}, is that $\left(\frac{-1}{n}\right) = (-1)^{(n-1)/2}$. Since we care about $n=2a-1 = 7 + 12k$ then 6$(-1)^{(2a-2)/2} = (-1)^{3+6k}=-1$,  showing that $-1$ is not a quadratic residue modulo $2a-1$. Therefore, in Case (ii) as well, no $\lambda$ can exist that solves \eqref{Equation with the pair of congruence relations}. 
\vskip2mm
This proves that $u_{nb}(K_1\# K_2)\ge 2$. On the other hand, $u_{nb}(K_1\#(-\bar K_2))=1$ as can be verified in Figure \ref{Connected sum of two two-bridge knots}. Taking $a=4$ gives the two-bridge knots $K_{7/2}=\bar 5_2$ and $K_{9/2}=\bar 6_1$. This completes the proof of Theorem \ref{Theorem about how the knot group does not determine the band unknotting number}.
\bibliographystyle{plain}
\bibliography{bibliography}

\end{document}